# A complete characterization on the robust isolated calmness of the nuclear norm regularized convex optimization problems

Ying Cui[*], Defeng Sun[†]

February 19, 2017

**Abstract**

In this paper, we provide a complete characterization on the robust isolated calmness of the Karush-Kuhn-Tucker (KKT) solution mapping for convex constrained optimization problems regularized by the nuclear norm function. This study is motivated by the recent work in [8], where the authors show that under the Robinson constraint qualification at a local optimal solution, the KKT solution mapping for a wide class of conic programming problem is robustly isolated calm if and only if both the second order sufficient condition (SOSC) and the strict Robinson constraint qualification (SRCQ) are satisfied. Based on the variational properties of the nuclear norm function and its conjugate, we establish the equivalence between the primal/dual SOSC and the dual/primal SRCQ. The derived results lead to several equivalent characterizations of the robust isolated calmness of the KKT solution mapping and add insights to the existing literature on the stability of the nuclear norm regularized convex optimization problems.

**Keywords.** robust isolated calmness, nuclear norm, second order sufficient condition, strict Robinson constraint qualification

**AMS subject classifications:** 90C25, 90C31, 65K10

## 1 Introduction

Let $\mathbb{X}$ and $\mathbb{Y}$ be two finite dimensional Euclidean spaces. Let $G : \mathbb{X} \rightrightarrows \mathbb{Y}$ be a set-valued mapping. The graph of $G$ is defined as $\operatorname{gph} G := \{(x, y) \in \mathbb{X} \times \mathbb{Y} \mid y \in G(x)\}$. Consider any $(\bar{x}, \bar{y}) \in \operatorname{gph} G$. The mapping $G$ is said to be isolated calm at $\bar{x}$ for $\bar{y}$ if there exist a constant $\kappa > 0$ and open neighborhoods $\mathcal{X}$ of $\bar{x}$ and $\mathcal{Y}$ of $\bar{y}$ such that

$$G(x) \cap \mathcal{Y} \subset \{\bar{y}\} + \kappa \|x - \bar{x}\| \mathbb{B}_{\mathbb{Y}}, \quad \forall\, x \in \mathcal{X}, \tag{1}$$

where $\mathbb{B}_{\mathbb{Y}}$ is the unit ball in $\mathbb{Y}$ (cf. e.g., [9, 3.9 (3I)]). The mapping $G$ is said to be robustly isolated calm at $\bar{x}$ for $\bar{y}$ if (1) holds and $G(x) \cap \mathcal{Y} \neq \emptyset$ for any $x \in \mathcal{X}$ [8, Definition 2].

[*]Department of Mathematics, National University of Singapore, 10 Lower Kent Ridge Road, Singapore (matcuiy@nus.edu.sg).

[†]Department of Mathematics, National University of Singapore, 10 Lower Kent Ridge Road, Singapore (matsundf@nus.edu.sg). This research is supported in part by the Academic Research Fund under Grant R-146-000-207-112.

In this paper, we are interested in characterizing the robust isolated calmness of the Karush-Kuhn-Tucker (KKT) solution mapping associated with the following nuclear norm regularized convex optimization problem:

$$\begin{array}{ll} \min\limits_{X} & h(\mathcal{F}X) + \langle C, X\rangle + \|X\|_* \\ \text{s.t.} & \mathcal{A}X - b \in \mathcal{Q}, \end{array} \tag{2}$$

where the function $h : \mathbb{R}^d \to \mathbb{R}$ is twice continuously differentiable on $\operatorname{dom} h$, which is assumed to be a non-empty open convex set, and is also essentially strictly convex (i.e., $h$ is strictly convex on every convex subset of $\operatorname{dom} \partial h$), $\mathcal{F} : \mathbb{R}^{m\times n} \to \mathbb{R}^d$ and $\mathcal{A} : \mathbb{R}^{m\times n} \to \mathbb{R}^e$ are linear operators, $C \in \mathbb{R}^{m\times n}$ and $b \in \mathbb{R}^e$ are given data, $\mathcal{Q} \subseteq \mathbb{R}^e$ is a nonempty convex polyhedral cone, $\|\cdot\|_*$ denotes the nuclear norm function in $\mathbb{R}^{m\times n}$, i.e., the sum of all the singular values of a given matrix, and $m, n, d, e$ are non-negative integers. The nuclear norm regularizer has been extensively used in diverse disciplines due to its ability in promoting a low rank solution. See the references [18, 19, 2, 3, 13, 16] for a sample of applications.

The concept of the isolated calmness is of fundamental importance in variational analysis. The monograph [9] by Dontchev and Rockafellar contains a comprehensive study on this subject. In the optimization field, the isolated calmness of the KKT solution mapping, besides its own interest in sensitivity analysis and perturbation theory, can be employed to investigate the convergence rates of primal dual type methods, including the proximal augmented Lagrangian method [15] and the alternating direction method of multipliers [10].

Obviously problem (2) can be equivalently formulated as the following conic programming problem

$$\begin{array}{ll} \min\limits_{X,t} & h(\mathcal{F}X) + \langle C, X\rangle + t \\ \text{s.t.} & \mathcal{A}X - b \in \mathcal{Q}, \quad (X, t) \in \operatorname{epi} \|\cdot\|_*, \end{array} \tag{3}$$

where $\operatorname{epi} \|\cdot\|_*$ denotes the epigraph of the function $\|\cdot\|_*$. Since $\operatorname{epi} \|\cdot\|_*$ is not a polyhedral set, the sensitivity results in the conventional nonlinear programming are not applicable for problem (3). Recently, some progress has been achieved in characterizing the isolated calmness of KKT solution mappings for problems involving non-polyhedral functions. For example, Zhang and Zhang [21] show that for the nonlinear semidefinite programming, the second order sufficient condition (SOSC) and the strict Robinson constraint qualification (SRCQ) at a local optimal solution together are sufficient for the KKT solution mapping to be isolated calm. Adding to this result, Han, Sun and Zhang [10] show that the SRCQ is also necessary to ensure the isolated calmness of the KKT solution mapping for such problems. In [11], Liu and Pan extend the aforementioned results to problems constrained by the epigraph of the Ky Fan $k$-norm function. The most recent work of Ding, Sun and Zhang [8] indicates that under the Robinson constraint qualification (RCQ) at a local optimal solution, the KKT solution mapping for a wide class of conic programming is robustly isolated calm at the origin for a KKT point if and only if both the SOSC and the SRCQ hold at the reference point.

The results developed in [8] can be directly applied to problem (3). Thus, by examining the relationships between the SOSCs, the (strict) RCQs as well as the robust isolated calmness of the solution mappings corresponding to problem (2) and problem (3), we are able to extend the work in [8] to the nuclear norm regularized convex optimization problem (2). Additionally, due to the special structure of problem (2) and its dual, we could provide more insightful characterizations about the isolated calmness of the KKT solution mapping. Note that the Lagrangian dual of problem (2) is

given by

$$\begin{array}{ll}\max\limits_{y,w,S} & -\langle b, y\rangle - \delta_{Q^\circ}(y) - h^*(w) \\ \text{s.t.} & \mathcal{A}^*y + \mathcal{F}^*w + S + C = 0, \quad \|S\|_2 \leqslant 1,\end{array} \tag{4}$$

where $\mathcal{A}^*$ and $\mathcal{F}^*$ are the adjoint of $\mathcal{A}$ and $\mathcal{F}$, respectively, $h^*(\cdot)$ is the conjugate function of $h$ and $\|\cdot\|_2$ denotes the spectral norm in $\mathbb{R}^{m\times n}$, i.e., the largest singular value of a given matrix. We shall show that for problem (2), its SRCQ is equivalent to the SOSC of problem (4), and conversely, its SOSC is equivalent to the SRCQ of problem (4). Armed with these results, we are led to a deep understanding of the robust isolated calmness for the KKT solution mapping for the nuclear norm regularized convex optimization problems.

The remaining parts of this paper are organized as follows. In the next section, we provide some preliminary results on variational analysis. In Section 3, we demonstrate how to translate the results of set-constrained problems in [8] into the language of nonsmooth optimization problems. In particular, this translation provides us a characterization of the robust isolated calmness of the KKT solution mapping for the nuclear norm regularized convex optimization problems. Section 4 is devoted to the study on the variational properties of the nuclear norm function. The derived results play an important role in our subsequent analysis. In Section 5, we establish the equivalence between the SOSC for the primal/dual problem and the SRCQ for the dual/primal problem. This establishment enables us to describe the robust isolated calmness of the KKT solution mapping for problem (2) via several equivalent conditions.

The following notation will be used throughout our paper.

- For a given positive integer $p$, we use $\mathbb{S}^p$ to denote the linear space of all $p \times p$ real symmetric matrices, $\mathbb{S}^p_+$ the cone of all $p \times p$ positive semidefinite matrices and $\mathbb{S}^p_-$ the cone of all $p \times p$ negative semidefinite matrices.

- For a given proper closed convex function $\theta : \mathbb{X} \to (-\infty, +\infty]$, we use $\operatorname{dom}\theta$ to denote its effective domain, $\operatorname{epi}\theta$ to denote its epigraph, $\theta^*$ to denote its conjugate, $\partial\theta$ to denote its subdifferential and $\operatorname{Prox}_\theta$ to denote its proximal mapping, all as in standard convex analysis [14].

- Let $D \subseteq \mathbb{R}^{m\times n}$ be a non-empty closed convex set. We write $\delta_D(\cdot)$ as the indicator function over $D$, i.e., $\delta_D(X) = 0$ if $X \in D$, and $\delta_D(X) = \infty$ if $X \notin D$. We write $\Pi_D(\cdot)$ as the metric projection onto $D$, i.e., $\Pi_D(X) := \arg\min_Y\{\|Y - X\| \mid Y \in D\}$ for $X \in \mathbb{R}^{m\times n}$.

- For any $z \in \mathbb{R}^m$, we denote $\operatorname{Diag}(z)$ as the $m \times m$ diagonal matrix whose $i$-th diagonal entry is $z_i$ for $i = 1, \ldots, m$. Let $\alpha \subseteq \{1, ..., m\}$ and $\beta \subseteq \{1, ..., n\}$ be two index sets. For any $Z \in \mathbb{R}^{m\times n}$, we write $Z_\alpha$ as the sub-matrix of $Z$ by removing all the columns of $Z$ not in $\alpha$, and $Z_{\alpha\beta}$ to be the $|\alpha| \times |\beta|$ sub-matrix of $Z$ obtained by removing all the rows of $Z$ not in $\alpha$ and all the columns of $Z$ not in $\beta$.

- Let $\mathcal{O}^n$ be the set of all $n \times n$ orthogonal matrices. For any $X \in \mathbb{R}^{m\times n}$, let $\sigma(X) \in \mathbb{R}^m$ be the vector of all singular values of $X$ with the entries $\sigma_1(X) \geqslant \sigma_2(X) \geqslant \ldots \geqslant \sigma_m(X)$, and let $\mathcal{O}^{m,n}(X)$ be the set of paired orthogonal matrices satisfying the singular value decomposition of $X$, i.e.,
$$\mathcal{O}^{m,n}(X) = \{(U, V) \in \mathcal{O}^m \times \mathcal{O}^n \mid X = U[\operatorname{Diag}(\sigma(X)) \;\; 0]V^T\}.$$

## 2 Preliminaries

In this section, we gather some knowledge on variational analysis that will be used in our subsequent developments. One can refer to the monograph [1] of Bonnans and Shapiro for detailed discussions on this subject.

A cone $\mathcal{Q} \subseteq \mathbb{Y}$ is said to be pointed if $y \in \mathcal{Q}$ and $-y \in \mathcal{Q}$ implies that $y = 0$. Let $\mathcal{Q} \subseteq \mathbb{Y}$ be a pointed convex closed cone. The closed convex set $\mathcal{K} \subseteq \mathbb{X}$ is said to be $\mathcal{C}^2$-cone reducible at $x \in \mathcal{K}$ to the cone $\mathcal{Q}$, if there exist an open neighborhood $\mathcal{W} \subseteq \mathbb{X}$ of $x$ and a twice continuously differentiable mapping $\Xi : \mathcal{W} \to \mathbb{Y}$ such that: (i) $\Xi(x) = 0 \in \mathbb{Y}$; (ii) the derivative mapping $\Xi(x) : \mathbb{X} \to \mathbb{Y}$ is onto; (iii) $\mathcal{K} \cap \mathcal{W} = \{x \in \mathcal{W} \mid \Xi(x) \in \mathcal{Q}\}$. We say that $\mathcal{K}$ is $\mathcal{C}^2$-cone reducible if $\mathcal{K}$ is $\mathcal{C}^2$-cone reducible at every $x \in \mathcal{K}$. A proper closed convex function $\theta : \mathbb{X} \to (-\infty, \infty]$ is said to be $\mathcal{C}^2$-cone reducible at $x \in \operatorname{dom}\theta$ if $\operatorname{epi}\theta$ is $\mathcal{C}^2$-cone reducible at $(x, \theta(x))$. Moreover, $\theta$ is said to be $\mathcal{C}^2$-cone reducible if it is $\mathcal{C}^2$-cone reducible at every $x \in \operatorname{dom}\theta$.

Given a subset $\mathcal{K} \subseteq \mathbb{X}$ and $x \in \mathcal{K}$, the contingent cone and the inner tangent cone of $\mathcal{K}$ at $x$ are defined as

$$\mathcal{T}_{\mathcal{K}}(x) = \limsup_{t \downarrow 0} \frac{\mathcal{K} - x}{t}$$

and

$$\mathcal{T}^i_{\mathcal{K}}(x) = \liminf_{t \downarrow 0} \frac{\mathcal{K} - x}{t},$$

respectively. If $\mathcal{K}$ is convex, $\dfrac{\mathcal{K} - x}{t}$ is a monotone decreasing function of $t$ such that $\mathcal{T}_{\mathcal{K}}(x) = \mathcal{T}^i_{\mathcal{K}}(x)$ for any $x \in \mathcal{K}$ [1, Proposition 2.55]. In this case, both $\mathcal{T}_{\mathcal{K}}(x)$ and $\mathcal{T}^i_{\mathcal{K}}(x)$ are called the tangent cone of $\mathcal{K}$ at $x$. Given $x \in \mathcal{K}$ and a direction $d \in \mathbb{X}$, define the inner and outer second order tangent sets at $x$ in the direction $d$ as

$$\mathcal{T}^{i,2}_{\mathcal{K}}(x; d) := \liminf_{t \downarrow 0} \frac{\mathcal{K} - x - td}{\frac{1}{2}t^2}$$

and

$$\mathcal{T}^{2}_{\mathcal{K}}(x; d) := \limsup_{t \downarrow 0} \frac{\mathcal{K} - x - td}{\frac{1}{2}t^2},$$

respectively. The inner and outer second order tangent sets not necessarily coincide in general, even if the set $\mathcal{K}$ is closed and convex. However, if $\mathcal{K}$ is a $\mathcal{C}^2$-cone reducible convex set, then $\mathcal{T}^{i,2}_{\mathcal{K}}(x; d) = \mathcal{T}^{2}_{\mathcal{K}}(x; d)$ [1, Proposition 3.136].

For a given function $\theta : \mathbb{X} \to (-\infty, +\infty]$, the lower and upper directional epiderivatives of $\theta$ at $x \in \operatorname{dom}\theta$ in the direction $h \in \mathbb{X}$ are defined as

$$\theta^{\downarrow}_{-}(x; h) := \liminf_{\substack{t \downarrow 0 \\ h' \to h}} \frac{\theta(x + th') - \theta(x)}{t}$$

and

$$\theta^{\downarrow}_{+}(x; h) := \sup_{\{t_n\} \in \Sigma} \left( \liminf_{\substack{n \to \infty \\ h' \to h}} \frac{\theta(x + t_n h') - \theta(x)}{t_n} \right),$$

respectively, where $\Sigma$ denotes the set of positive real sequences $\{t_n\}$ converging to 0. The contingent and inner tangent cone of $\operatorname{epi}\theta$ are closely related to the lower and upper directional epiderivative of $\theta$ [1, proposition 2.58]. Specifically, for any $x \in \operatorname{dom}\theta$,

$$\mathcal{T}_{\operatorname{epi}\theta}\big(x, \theta(x)\big) = \operatorname{epi}\theta^{\downarrow}_{-}(x; \cdot), \quad \mathcal{T}^i_{\operatorname{epi}\theta}\big(x, \theta(x)\big) = \operatorname{epi}\theta^{\downarrow}_{+}(x; \cdot). \tag{5}$$

One can observe from the above equations that if $\theta$ is a closed convex function, then $\theta^{\downarrow}_{-}(x;\cdot) = \theta^{\downarrow}_{+}(x;\cdot)$ for any $x \in \operatorname{dom}\theta$. In this case we say that $\theta$ is directionally epidifferentiable at $x$ and write the common value as $\theta^{\downarrow}(x;\cdot)$. If $\theta^{\downarrow}_{+}(x;d)$ and $\theta^{\downarrow}_{-}(x;d)$ are finite for $x \in \operatorname{dom}\theta$ and $d \in \mathbb{X}$, we also define the following lower and upper second order epiderivatives for $w \in \mathbb{X}$:

$$\theta^{\downarrow\downarrow}_{-}(x;d,w) := \liminf_{\substack{t\downarrow 0\\ w'\to w}} \frac{\theta(x+td+\frac{1}{2}t^2w') - \theta(x) - t\theta^{\downarrow}_{-}(x;d)}{\frac{1}{2}t^2},$$
$$\theta^{\downarrow\downarrow}_{+}(x;d,w) := \sup_{t_n\in\Sigma}\left(\liminf_{\substack{n\to\infty\\ w'\to w}} \frac{\theta(x+t_nd+\frac{1}{2}t_n^2w') - \theta(x) - t_n\theta^{\downarrow}_{+}(x;d)}{\frac{1}{2}t_n^2}\right).$$

Similarly to (5), the inner and outer second order tangent sets of epi $\theta$ are closely related to the lower and upper second order epiderivative of $\theta$ [1, proposition 3.41]. Specifically, for any $x \in \operatorname{dom}\theta$ and $d \in \mathbb{X}$, if $\theta^{\downarrow}_{+}(x;d)$ and $\theta^{\downarrow}_{-}(x;d)$ are finite, then

$$\mathcal{T}^{i,2}_{\operatorname{epi}\theta}((x,\theta(x));(d,\theta^{\downarrow}_{+}(x;d)) = \operatorname{epi}\theta^{\downarrow\downarrow}_{+}(x;d,\cdot), \quad \mathcal{T}^{2}_{\operatorname{epi}\theta}((x,\theta(x));(d,\theta^{\downarrow}_{-}(x;d)) = \operatorname{epi}\theta^{\downarrow\downarrow}_{-}(x;d,\cdot). \tag{6}$$

# 3 Constraint qualifications, second order sufficient optimality conditions and robust isolated calmness

Consider the following canonical perturbation of a general class of nonsmooth optimization problems (not necessarily convex):

$$\min\left\{f(x) + \theta(x) - \langle \delta_1, x\rangle \mid g(x) + \delta_2 \in \mathcal{P}\right\}, \tag{7}$$

where $f : \mathbb{X} \to \mathbb{R}$ and $g : \mathbb{X} \to \mathbb{Y}$ are twice continuously differentiable functions, $\mathcal{P} \subseteq \mathbb{Y}$ is a closed convex set, $\theta : \mathbb{X} \to (-\infty, +\infty]$ is a closed proper convex function, and $\delta_1 \in \mathbb{X}$ and $\delta_2 \in \mathbb{Y}$ are perturbation parameters.

Note that problem (7) can be equivalently written as the following optimization problem:

$$\min\left\{f(x) + t - \langle \delta_1, x\rangle \mid g(x) + \delta_2 \in \mathcal{P},\ (x,t) \in \mathcal{K}\right\}, \tag{8}$$

where $\mathcal{K} := \operatorname{epi}\theta$ is a closed convex set. The constraint qualifications and SOSCs for problem (8) have been extensively explored in [1, Section 3], through the study of the (second order) tangent sets of $\mathcal{P}$ and $\mathcal{K}$ at a stationary point. In the following, by employing the equations in (5) and (6), we reduce these properties to the (second order) directional epiderivatives of $\theta$. This reduction leads to a direct approach to the sensitivity analysis of the nonsmooth optimization problem (7).

For notational simplicity, denote $\mathbb{Z} := \mathbb{X} \times \mathbb{R} \times \mathbb{Y} \times \mathbb{X} \times \mathbb{R}$. Let $(\delta_1, \delta_2) \in \mathbb{X} \times \mathbb{Y}$ be given. We say that $(\bar{x}, \bar{t})$ is a feasible solution to problem (8) if

$$(\bar{x}, \bar{t}) \in \widehat{F}(\delta_1, \delta_2) := \{(x,t) \in \mathbb{X} \times \mathbb{R} \mid\ g(x) + \delta_2 \in \mathcal{P},\ (x,t) \in \mathcal{K}\}.$$

For any $(x,t,y,z,\tau) \in \mathbb{Z}$, the Lagrangian function of (8) with $(\delta_1,\delta_2) = 0$ can be written as

$$\mathcal{L}(x,t;y,z,\tau) := f(x) + t + \langle y, g(x)\rangle + \langle z, x\rangle + t\tau.$$

For any given $(\delta_1, \delta_2) \in \mathbb{X} \times \mathbb{Y}$, the KKT optimality condition for problem (8) is

$$\begin{cases} \delta_1 = \nabla_x \mathcal{L}(x,t,y,z,\tau) = \nabla f(x) + \nabla g(x)y + z, \\ \nabla_t \mathcal{L}(x,t,y,z,\tau) = 1 + \tau = 0, \qquad\qquad (x,t,y,z,\tau) \in \mathbb{Z}, \\ y \in \mathcal{N}_{\mathcal{P}}(g(x) + \delta_2), \quad (z,\tau) \in \mathcal{N}_{\mathcal{K}}(x,t), \end{cases} \tag{9}$$

where $\mathcal{N}_C(s)$ denotes the normal cone of a given convex set $C$ at $s \in C$. Let $\widehat{S}_{\rm KKT} : \mathbb{X} \times \mathbb{Y} \to \mathbb{Z}$ be the following KKT solution mapping:

$$\widehat{S}_{\rm KKT}(\delta_1, \delta_2) := \{(x, t, y, z, \tau) \in \mathbb{Z} \mid (x, t, y, z, \tau) \text{ satisfies (9)}\}, \quad (\delta_1, \delta_2) \in \mathbb{X} \times \mathbb{Y}. \tag{10}$$

For any given $(\delta_1, \delta_2) \in \mathbb{X} \times \mathbb{Y}$, we call $(\bar{x}, \bar{t})$ a stationary point of problem (7) if there exists a Lagrangian multiplier $(\bar{y}, \bar{z}, \bar{\tau}) \in \mathbb{Y} \times \mathbb{X} \times \mathbb{R}$ such that $(\bar{x}, \bar{t}, \bar{y}, \bar{z}, \bar{\tau}) \in \widehat{S}_{\rm KKT}(\delta_1, \delta_2)$. Denote $\widehat{\mathcal{M}}(\bar{x}, \bar{t}, \delta_1, \delta_2) \subseteq \mathbb{Y} \times \mathbb{X} \times \mathbb{R}$ as the set of all such Lagrangian multipliers $(\bar{y}, \bar{z}, \bar{\tau})$ associated with $(\bar{x}, \bar{t})$. Note from (9) that $\bar{\tau} \equiv -1$.

Let $(\delta_1, \delta_2) = 0$ in problem (8). The RCQ is said to hold at a feasible solution $(\bar{x}, \theta(\bar{x}))$ of problem (8) if

$$\begin{pmatrix} (g'(\bar{x}), 0) \\ (\mathcal{I}_{\mathbb{X}}, 1) \end{pmatrix} (\mathbb{X} \times \mathbb{R}) + \begin{pmatrix} \mathcal{T}_{\mathcal{P}}(g(\bar{x})) \\ \mathcal{T}_{\mathcal{K}}(\bar{x}, \theta(\bar{x})) \end{pmatrix} = \begin{pmatrix} \mathbb{Y} \\ \mathbb{X} \times \mathbb{R} \end{pmatrix}, \tag{11}$$

where $\mathcal{I}_{\mathbb{X}}$ is the identity mapping from $\mathbb{X}$ to $\mathbb{X}$. It is known that the RCQ (11) holds at a local optimal solution $(\bar{x}, \theta(\bar{x}))$ if and only if $\widehat{\mathcal{M}}(\bar{x}, \theta(\bar{x}), 0, 0)$ is a nonempty, convex and compact set (cf., e.g., [1, Theorem 3.9 and Proposition 3.17]). The SRCQ is said to hold at $(\bar{x}, \theta(\bar{x}))$ for $(\bar{y}, \bar{z}, -1) \in \widehat{\mathcal{M}}(\bar{x}, \theta(\bar{x}), 0, 0)$ if

$$\begin{pmatrix} (g'(\bar{x}), 0) \\ (\mathcal{I}_{\mathbb{X}}, 1) \end{pmatrix} (\mathbb{X} \times \mathbb{R}) + \begin{pmatrix} \mathcal{T}_{\mathcal{P}}(g(\bar{x})) \cap \bar{y}^{\perp} \\ \mathcal{T}_{\mathcal{K}}(\bar{x}, \theta(\bar{x})) \cap (\bar{z}, -1)^{\perp} \end{pmatrix} = \begin{pmatrix} \mathbb{Y} \\ \mathbb{X} \times \mathbb{R} \end{pmatrix}. \tag{12}$$

Obviously the SRCQ (12) is stronger than the RCQ (11). It follows from [1, Proposition 4.50] that $\widehat{M}(\bar{x}, \theta(\bar{x}), 0, 0)$ is a singleton if the SRCQ holds at $(\bar{x}, \theta(\bar{x}))$. The critical cone at a feasible point $(\bar{x}, \theta(\bar{x}))$ for problem (8) takes the form of

$$\widehat{\mathcal{C}}(\bar{x}, \theta(\bar{x})) := \{(d_1, d_2) \in \mathbb{X} \times \mathbb{R} \mid g'(\bar{x})d_1 \in \mathcal{T}_{\mathcal{P}}(g(\bar{x})),\ (d_1, d_2) \in \mathcal{T}_{\mathcal{K}}\big(\bar{x}, \theta(\bar{x})\big),\ f'(\bar{x})d_1 + d_2 \leqslant 0\}.$$

Furthermore, if $(\bar{x}, \theta(\bar{x}))$ is a stationary point of problem (8) and there exists $(\bar{y}, \bar{z}, -1) \in \widehat{\mathcal{M}}(\bar{x}, \theta(\bar{x}), 0, 0)$, then

$$\begin{aligned} \widehat{\mathcal{C}}(\bar{x}, \theta(\bar{x})) &= \{(d_1, d_2) \in \mathbb{X} \times \mathbb{R} \mid g'(\bar{x})d_1 \in \mathcal{T}_{\mathcal{P}}(g(\bar{x})),\ (d_1, d_2) \in \mathcal{T}_{\mathcal{K}}\big(\bar{x}, \theta(\bar{x})\big),\ f'(\bar{x})d_1 + d_2 = 0\} \\ &= \{(d_1, d_2) \in \mathbb{X} \times \mathbb{R} \mid g'(\bar{x})d_1 \in \mathcal{T}_{\mathcal{P}}(g(\bar{x})) \cap \bar{y}^{\perp},\ (d_1, d_2) \in \mathcal{T}_{\mathcal{K}}\big(\bar{x}, \theta(\bar{x})\big) \cap (\bar{z}, -1)^{\perp}\}. \end{aligned}$$

Assume that $\operatorname{epi}\theta$ is second order regular at $(\bar{x}, \theta(\bar{x}))$ (see [1, Definition 3.85] for the definition of the second order regularity). This assumption is particularly satisfied when $\theta(\cdot) = \|\cdot\|_*$, since in this case $\operatorname{epi}\theta$ is a $\mathcal{C}^2$-cone reducible set [5, Proposition 4.3]. Then the SOSC at a stationary point $(\bar{x}, \theta(\bar{x}))$ for problem (8) with $(\delta_1, \delta_2) = 0$ is said to hold if for any $d \in \widehat{\mathcal{C}}(\bar{x}, \theta(\bar{x})) \backslash \{0\}$,

$$\sup_{(\bar{y}, \bar{z}, -1) \in \widehat{\mathcal{M}}(\bar{x}, \theta(\bar{x}), 0, 0)} \langle d, \nabla^2_{(x,t)(x,t)} \mathcal{L}(\bar{x}, \theta(\bar{x}), \bar{y}, \bar{z}, -1) d\rangle - \sigma\big((\bar{z}, -1), \mathcal{T}^2_{\mathcal{K}}((\bar{x}, \theta(\bar{x})); (g'(\bar{x}), 0)d)\big) > 0, \tag{13}$$

where $\sigma(s, C) := \sup\{\langle s', s\rangle \mid s' \in C\}$ denotes the support function of a given set $C$ at $s$. The above SOSC implies the quadratic growth condition at $(\bar{x}, \theta(\bar{x}))$ (cf. e.g., [1, Theorem 3.86]), that is, there exist a constant $\kappa > 0$ and a neighborhood $N$ of $(\bar{x}, \theta(\bar{x}))$ such that

$$f(x) + t \geqslant f(\bar{x}) + \theta(\bar{x}) + \kappa \, \|(x, t) - (\bar{x}, \theta(\bar{x}))\|^2, \quad \forall\ (x, t) \in \widehat{F}(0, 0) \cap N.$$

The following proposition, which is taken from [8, Theorem 24], characterizes the robust isolated calmness of problem (8) via the SOSC (13) and the SRCQ (12).

**Proposition 3.1.** *Suppose that $(\bar{x}, \theta(\bar{x})) \in \mathbb{X}\times\mathbb{R}$ is a feasible solution of problem (8) with $(\delta_1, \delta_2) = 0$. Suppose that the RCQ (11) holds at $(\bar{x}, \theta(\bar{x}))$. Assume that $\operatorname{epi}\theta$ is $\mathcal{C}^2$-cone reducible at $(\bar{x}, \theta(\bar{x}))$. Let $(\bar{y}, \bar{z}, -1) \in \widehat{\mathcal{M}}(\bar{x}, \theta(\bar{x}), 0, 0)$. Then the following two statements are equivalent to each other:*
*(i) The SOSC (13) holds at $(\bar{x}, \theta(\bar{x}))$ and the SRCQ (12) holds at $(\bar{x}, \theta(\bar{x}))$ for $(\bar{y}, \bar{z}, -1)$.*
*(ii) The point $(\bar{x}, \theta(\bar{x}))$ is a local optimal solution of problem (8) and the KKT solution mapping $\widehat{S}_{KKT}$ is robustly isolated calm at the origin for $(\bar{x}, \theta(\bar{x}), \bar{y}, \bar{z}, -1)$.*

Now we return to the nonsmooth optimization problem (7). Let $(\delta_1, \delta_2) \in \mathbb{X} \times \mathbb{Y}$ be given. We say that $\bar{x}$ is a feasible solution to problem (7) if

$$\bar{x} \in F(\delta_1, \delta_2) := \{x \in \operatorname{dom}\theta \mid \; g(x) + \delta_2 \in \mathcal{P}\}.$$

Denote $l : \mathbb{X} \times \mathbb{Y} \to \mathbb{R}$ by

$$l(x, y) := f(x) + \langle g(x), y\rangle, \quad (x, y) \in \mathbb{X} \times \mathbb{Y}.$$

Then the KKT optimality condition takes the form of

$$\begin{cases} \delta_1 \in \nabla_x l(x, y) + \partial\theta(x), \\ y \in \mathcal{N}_{\mathcal{P}}(g(x) + \delta_2), \end{cases} \quad (x, y) \in \mathbb{X} \times \mathbb{Y}. \tag{14}$$

Let $S_{\rm KKT} : \mathbb{X} \times \mathbb{Y} \to \mathbb{X} \times \mathbb{Y}$ be the following KKT solution mapping:

$$S_{\rm KKT}(\delta_1, \delta_2) := \{(x, y) \in \mathbb{X} \times \mathbb{Y} \mid (x, y) \text{ satisfies (14)}\},\ (\delta_1, \delta_2) \in \mathbb{X} \times \mathbb{Y}. \tag{15}$$

For any given $(\delta_1, \delta_2) \in \mathbb{X}\times\mathbb{Y}$, we call $\bar{x}$ a stationary point of problem (7) if there exists a Lagrangian multiplier $\bar{y} \in \mathbb{Y}$ such that $(\bar{x}, \bar{y}) \in S_{\rm KKT}(\delta_1, \delta_2)$. Denote $\mathcal{M}(\bar{x}, \delta_1, \delta_2) \subseteq \mathbb{Y}$ as the set of all such Lagrangian multipliers $\bar{y}$ associated with $\bar{x}$.

The following proposition establishes the equivalence between the robust isolated calmness of the KKT solution mappings with respect to problem (7) and problem (8).

**Proposition 3.2.** *Let $(\bar{x}, \theta(\bar{x})) \in \mathbb{X}\times\mathbb{R}$ be a local optimal solution of problem (8) with $\widehat{\mathcal{M}}(\bar{x}, \theta(\bar{x}), 0, 0) \neq \emptyset$. Let $(\bar{y}, \bar{z}, -1) \in \widehat{M}(\bar{x}, \theta(\bar{x}), 0, 0)$. If the KKT solution mapping $\widehat{S}_{KKT}$ given in (10) is robustly isolated calm at the origin for $(\bar{x}, \theta(\bar{x}), \bar{y}, \bar{z}, -1)$, then the KKT solution mapping $S_{KKT}$ given in (15) is robustly isolated calm at the origin for $(\bar{x}, \bar{y})$. The reverse implication is true if the function $\theta$ is Lipschitz continuous at $\bar{x}$.*

*Proof.* Note from [4, Corollary 2.4.9] that

$$(z, -1) \in \mathcal{N}_{\operatorname{epi}\theta}\big(x, \theta(x)\big) \Longleftrightarrow z \in \partial\theta(x), \quad \forall\, x, z \in \mathbb{X}.$$

Then for any $(\delta_1, \delta_2) \in \mathbb{X} \times \mathbb{Y}$ and any $(x, t, y, z, -1) \in \widehat{S}_{\rm KKT}(\delta_1, \delta_2)$, we know from (9) and (14) that $(x, y) \in S_{\rm KKT}(\delta_1, \delta_2)$. Thus, the first part of this proposition follows easily from the definition of the robust isolated calmness.

Conversely, consider any $(\delta_1, \delta_2) \in \mathbb{X}\times\mathbb{Y}$ and $(x, y) \in S_{\rm KKT}(\delta_1, \delta_2)$. Let $z = \delta_1 - \nabla f(x) - \nabla g(x)y$. By the similar arguments as above, we have $(x, \theta(x), y, z, -1) \in \widehat{S}_{\rm KKT}(\delta_1, \delta_2)$. Since $f$ and $g$ are assumed to be twice continuously differentiable, $\nabla f(\cdot)$ and $\nabla g(\cdot)$ are locally Lipschitz continuous at $\bar{x}$. Then there exists a constant $\kappa > 0$ (only depending on $\bar{x}$ and $\bar{y}$) such that for any $x$ sufficiently close to $\bar{x}$,

$$\begin{cases} \|\theta(x) - \theta(\bar{x})\| \leqslant \kappa\|x - \bar{x}\|, \\ \|z - \bar{z}\| \leqslant \|\delta_1\| + \|\nabla f(x) - \nabla f(\bar{x})\| + \|\nabla g(x)y - \nabla g(\bar{x})\bar{y}\| \leqslant k(\|\delta_1\| + \|x - \bar{x}\| + \|y - \bar{y}\|). \end{cases}$$

Consequently, the second assertion of this proposition also follows from the definition of the robust isolated calmness. □

As mentioned in Section 2, the closed convex function $\theta(\cdot)$ is always directional epidifferentiable at $x \in \operatorname{dom}\theta$. Then by [14, Theorem 23.2], the KKT optimality condition (14) is equivalent to

$$\begin{cases} \theta^{\downarrow}(x; d) + \langle \nabla_x l(x, y) - \delta_1, d\rangle \geqslant 0, \\ y \in \mathcal{N}_{\mathcal{P}}(g(x) + \delta_2), \end{cases} \qquad \forall\, d \in \mathbb{X},$$

where $(x, y) \in \mathbb{X} \times \mathbb{Y}$. Define the critical of the function $\theta$ by

$$\mathcal{C}_\theta(x, z) := \{d \in \mathcal{T}_{\operatorname{dom}\theta}(x) \mid \theta^{\downarrow}(x; d) = \langle d, z\rangle\}, \quad (x, z) \in \operatorname{dom}\theta \times \mathbb{X}. \tag{16}$$

Let $(\delta_1, \delta_2) = 0$. The RCQ is said to hold at a feasible solution $\bar{x}$ of problem (7) if

$$\begin{pmatrix} g'(\bar{x}) \\ \mathcal{I}_{\mathbb{X}} \end{pmatrix} \mathbb{X} + \begin{pmatrix} \mathcal{T}_{\mathcal{P}}(g(\bar{x})) \\ \mathcal{T}_{\operatorname{dom}\theta}(\bar{x}) \end{pmatrix} = \begin{pmatrix} \mathbb{Y} \\ \mathbb{X} \end{pmatrix}. \tag{17}$$

By the equations in (5), the SRCQ is said to hold at a stationary point $\bar{x}$ for $\bar{y} \in \mathcal{M}(\bar{x}, 0, 0)$ if

$$\begin{pmatrix} g'(\bar{x}) \\ \mathcal{I}_{\mathbb{X}} \end{pmatrix} \mathbb{X} + \begin{pmatrix} \mathcal{T}_{\mathcal{P}}(g(\bar{x})) \cap \bar{y}^{\perp} \\ \mathcal{C}_\theta(\bar{x}, -\nabla_x l(\bar{x}, \bar{y})) \end{pmatrix} = \begin{pmatrix} \mathbb{Y} \\ \mathbb{X} \end{pmatrix}. \tag{18}$$

The critical cone at a feasible solution $\bar{x}$ of problem (7) is given by

$$\mathcal{C}(\bar{x}) := \left\{d \in \mathbb{X} \mid g'(\bar{x})d \in \mathcal{T}_{\mathcal{P}}(g(\bar{x})),\, d \in \mathcal{T}_{\operatorname{dom}\theta}(\bar{x}),\, f'(\bar{x})d + \theta^{\downarrow}(\bar{x}; d) \leqslant 0\right\}.$$

If $\bar{x}$ is a stationary point of problem (7) and $\bar{y} \in \mathcal{M}(\bar{x}, 0, 0)$, then

$$\begin{aligned} \mathcal{C}(\bar{x}) &= \left\{d \in \mathbb{X} \mid g'(\bar{x})d \in \mathcal{T}_{\mathcal{P}}(g(\bar{x})),\, d \in \mathcal{T}_{\operatorname{dom}\theta}(\bar{x}),\, f'(\bar{x})d + \theta^{\downarrow}(\bar{x}; d) = 0\right\} \\ &= \left\{d \in \mathbb{X} \mid g'(\bar{x})d \in \mathcal{T}_{\mathcal{P}}(g(\bar{x})) \cap \bar{y}^{\perp},\, d \in \mathcal{C}_\theta(\bar{x}, -\nabla f(\bar{x}))\right\}. \end{aligned}$$

Based on the equations in (6), we know that if $\theta$ is $\mathcal{C}^2$-cone reducible at $\bar{x}$, then the SOSC at $\bar{x}$ for problem (7) with $(\delta_1, \delta_2) = 0$ takes the form of

$$\sup_{\bar{y} \in \mathcal{M}(\bar{x}, 0, 0)} \left\{\langle d, \nabla^2_{xx} l(\bar{x}, \bar{y})d\rangle - \psi^*_{(\bar{x}, d)}(-\nabla_x l(\bar{x}, \bar{y}))\right\} > 0, \quad \forall\, d \in \mathcal{C}(\bar{x})\backslash\{0\}, \tag{19}$$

where $\psi^*_{(x,d)}(\cdot)$ is the conjugate function of $\psi_{(x,d)}(\cdot) = \theta^{\downarrow\downarrow}_{-}(x; d, \cdot)$ for any $x \in \operatorname{dom}\theta$ and any $d \in \mathbb{X}$.

Combining Proposition 3.1 and Proposition 3.2, we are ready to state the main result of this section.

**Proposition 3.3.** *Let $\bar{x} \in \mathbb{X}$ be a local optimal solution of problem (7) with $(\delta_1, \delta_2) = 0$. Suppose that the RCQ (17) holds at $\bar{x}$. Assume that $\theta$ is $\mathcal{C}^2$-cone reducible and Lipschitz continuous at $\bar{x}$. Let $\bar{y} \in \mathcal{M}(\bar{x}, 0, 0)$. Then the following two statements are equivalent to each other:*
*(i) The SOSC (19) holds at $\bar{x}$ and the SRCQ (18) holds at $\bar{x}$ for $\bar{y}$.*
*(ii) The point $\bar{x}$ is a local optimal solution of problem (7) and the KKT solution mapping $S_{KKT}$ is robustly isolated calm at the origin for $(\bar{x}, \bar{y})$.*

## 4 Variational analysis of the nuclear norm function

Throughout this section, we denote $\theta : \mathbb{R}^{m\times n} \to \mathbb{R}$ as the nuclear norm function. Obviously $\operatorname{dom}\theta = \mathbb{R}^{m\times n}$. Since the nuclear norm function is convex and globally Lipschitz continuous, for any $X \in \mathbb{R}^{m\times n}$, both $\theta_-^{\downarrow}(X;\cdot)$ and $\theta_+^{\downarrow}(X;\cdot)$ defined in Section 2 coincide with $\theta'(X;\cdot)$, the conventional directional derivative of $\theta$ at $X$ [1, Section 2.2.3].

Let $X \in \mathbb{R}^{m\times n}$ be an arbitrary but fixed point. Suppose that $X$ admits the following singular-value decomposition (SVD):

$$X = U[\operatorname{Diag}(\sigma(X))\ 0]V^T = U\operatorname{Diag}(\sigma(X))V_1^T, \tag{20}$$

where $U \in \mathcal{O}^m$ and $V = [V_1\ V_2] \in \mathcal{O}^n$ are the left and right singular vectors of $X$ with $V_1 \in \mathbb{R}^{n\times m}$ and $V_2 \in \mathbb{R}^{n\times(n-m)}$. Define the index sets

$$a := \{1 \leqslant i \leqslant m \mid \sigma_i(X) \in (1,+\infty)\}, \quad b := \{1 \leqslant i \leqslant m \mid \sigma_i(X) \in [0,1]\}, \quad c = \{m+1,\ldots,n\}. \tag{21}$$

Denote the distinct singular values of $X$ that are greater than 1 by $\nu_1(X) > \ldots > \nu_r(X) > 1$, where $r$ is a non-negative integer. We further divide the sets $a$ and $b$ into the following subsets:

$$\begin{aligned} &a_l := \{i \in a \mid \sigma_i(X) = \nu_l(X)\}, \quad l = 1,\ldots,r, \\ &b_1 := \{i \in b \mid \sigma_i(X) = 1\}, \quad b_2 := \{i \in b \mid 0 < \sigma_i(X) < 1\}, \quad b_3 := \{i \in b \mid \sigma_i(X) = 0\}. \end{aligned} \tag{22}$$

Let us first review the concept of Löwner's operator and its differential properties. Suppose that $X \in \mathbb{R}^{m\times n}$ has the SVD (20). For any scalar function $g : \mathbb{R} \to \mathbb{R}$, define the corresponding matrix valued function $G$ by

$$G(X) := U[\operatorname{Diag}(g(\sigma_1(X)), g(\sigma_2(X)),\ldots, g(\sigma_m(X))\ 0]V^T.$$

Such a function is called Löwner's operator associated with the function $g$, which is first studied by Löwner in the context of symmetric matrices [12]. In particular, let $\phi : \mathbb{R} \to \mathbb{R}$ be the scalar function

$$\phi(x) := \max\{x-1, 0\}, \quad x \in \mathbb{R}.$$

It is easy to verify that the proximal mapping of $\theta$ can be expressed as:

$$\operatorname{Prox}_\theta(X) = U[\operatorname{Diag}(\phi(\sigma_1(X)),\ldots,\phi(\sigma_m(X)))\ 0]V^T, \quad X \in \mathbb{R}^{m\times n}. \tag{23}$$

Clearly $\operatorname{Prox}_\theta(\cdot)$ can be taken as Löwner's operator associated with the function $\phi$. The directional derivate of $\operatorname{Prox}_\theta(\cdot)$ can thus be obtained via the general formula regarding the directional derivative of Löwner's operator [7]. Obviously $\phi$ is directionally differentiable with the directional derivative

$$\phi'(x;d) = \begin{cases} d & \text{if } x > 1, \\ \max\{d,0\} & \text{if } x = 1, \\ 0 & \text{if } x < 1, \end{cases} \qquad x \in \mathbb{R},\ d \in \mathbb{R}.$$

For any positive integer $p$, define linear matrix operators $S : \mathbb{R}^{p\times p} \to \mathbb{S}^p$ and $T : \mathbb{R}^{p\times p} \to \mathbb{R}^{p\times p}$ by

$$S(X) := \frac{1}{2}(X + X^T), \quad T(X) := \frac{1}{2}(X - X^T), \quad X \in \mathbb{R}^{p\times p}. \tag{24}$$

Denote $\Xi^2_{aa} : \mathbb{R}^{m\times n} \to \mathbb{R}^{|a|\times|a|}$, $\Xi^1_{ab} : \mathbb{R}^{m\times n} \to \mathbb{R}^{|a|\times|b|}$, $\Xi^2_{ab} : \mathbb{R}^{m\times n} \to \mathbb{R}^{|a|\times|b|}$ and $\Xi_2 : \mathbb{R}^{m\times n} \to \mathbb{R}^{|a|\times|c|}$ as

$$\begin{cases} ((\Xi^2_{aa})(X))_{ij} := \dfrac{\sigma_i(X)+\sigma_j(X)-2}{\sigma_i(X)+\sigma_j(X)}, \ i=1,2,\ldots,|a|, \ j=1,2,\ldots,|a|, \\ ((\Xi^1_{ab})(X))_{ij} := \dfrac{\sigma_i(X)-1}{\sigma_i(X)-\sigma_{j+|a|}(X)}, \ i=1,2,\ldots,|a|, \ j=1,2,\ldots,|b|, \\ ((\Xi^2_{ab})(X))_{ij} := \dfrac{\sigma_i(X)-1}{\sigma_i(X)+\sigma_{j+|a|}(X)}, \ i=1,2,\ldots,|a|, \ j=1,2,\ldots,|b|, \\ ((\Xi_2)(X))_{ij} := \dfrac{\sigma_i(X)-1}{\sigma_i(X)}, \quad i=1,2,\ldots,|a|, \ j=1,2,\ldots,n-m, \end{cases} \qquad X \in \mathbb{R}^{m\times n}.$$

Denote $\Gamma_1 : \mathbb{R}^{m\times n}\times\mathbb{R}^{m\times n} \to \mathbb{R}^{|a|\times|a|}$, $\Gamma_2 : \mathbb{R}^{m\times n}\times\mathbb{R}^{m\times n} \to \mathbb{R}^{|a|\times|b|}$, $\Gamma_3 : \mathbb{R}^{m\times n}\times\mathbb{R}^{m\times n} \to \mathbb{R}^{|b|\times|a|}$ and $\Gamma_4 : \mathbb{R}^{m\times n}\times\mathbb{R}^{m\times n} \to \mathbb{R}^{|a|\times|c|}$ as

$$\begin{cases} \Gamma_1(X,H) &:= & (S(H_1))_{aa} + \Xi^2_{aa}(X)\circ(T(H_1))_{aa}, \\ \Gamma_2(X,H) &:= & \Xi^1_{ab}(X)\circ(S(H_1))_{ab} + \Xi^2_{ab}(X)\circ(T(H_1))_{ab}, \\ \Gamma_3(X,H) &:= & (\Xi^1_{ab}(X))^T\circ(S(H_1))_{ba} + (\Xi^2_{ab}(X))^T\circ(T(H_1))_{ba}, \\ \Gamma_4(X,H) &:= & \Xi_2(X)\circ H_{ac}, \end{cases} \qquad (X,H)\in\mathbb{R}^{m\times n}\times\mathbb{R}^{m\times n},$$

where $\circ$ denotes the Hadamard product between two matrices and $H = [H_1 \ \ H_2]$ with $H_1 \in \mathbb{R}^{m\times m}$ and $H_2 \in \mathbb{R}^{m\times(n-m)}$. Then by [7, Theorem 3], the directional derivative of $\mathrm{Prox}_\theta(\cdot)$ at $X \in \mathbb{R}^{m\times n}$ in the direction $H \in \mathbb{R}^{m\times n}$ takes the form of

$$\mathrm{Prox}'_\theta(X;H) = U \left( \begin{array}{c|ccc|c} \Gamma_1(X,\widetilde{H}) & & \Gamma_2(X,\widetilde{H}) & & \Gamma_4(X,\widetilde{H}) \\ \hline & \Pi_{\mathbb{S}^{|b_1|}_+}(S(\widetilde{H}_{b_1b_1})) & 0 & 0 & \\ \Gamma_3(X,\widetilde{H}) & 0 & 0_{b_2\times b_2} & 0 & 0 \\ & 0 & 0 & 0_{b_3\times b_3} & \end{array} \right) V^T, \tag{25}$$

where $\widetilde{H} = U^T H V$.

In [17], Watson shows that the subdifferential of the nuclear norm function takes the following form:

$$\partial\theta(X) = \left\{ U_a V_a^T + U_b W [V_b \ V_2]^T \mid W \in \mathbb{R}^{|b|\times(n-|a|)}, \ \|W\|_2 \leqslant 1 \right\}, \quad X \in \mathbb{R}^{m\times n}. \tag{26}$$

Therefore, for any $H \in \mathbb{R}^{m\times n}$, the directional derivative of $\theta$ at $X$ in the direction $H$ is given by

$$\theta'(X;H) = \sup_{S\in\partial\theta(X)} \langle H, S\rangle = \mathrm{tr}(U_a^T H V_a) + \|U_b^T H [V_b \ V_2]\|_*. \tag{27}$$

Let $A = \mathrm{Prox}_\theta(X)$ and $B = \mathrm{Prox}_{\theta^*}(X)$. Define the critical cone of $\theta$ at $A$ for $B$ as

$$\mathcal{C}_\theta(A,B) := \{H \in \mathbb{R}^{m\times n} \mid \theta'(A;H) = \langle H, B\rangle\}. \tag{28}$$

Similarly, define the critical cone of $\theta^*$ at $B$ for $A$ as

$$\mathcal{C}_{\theta^*}(B,A) := \{H \in \mathbb{R}^{m\times n} \mid (\theta^*)'(B;H) = \langle H, A\rangle\}. \tag{29}$$

As can be seen from (19), in order to analyze the SOSC for problem (2), one needs to compute the conjugate of the second order epiderivative of $\theta$. This has already been done in [6]. Firstly, it follows from (23) that $\sigma(A) = \max\{\sigma(X) - 1, 0\}$. Specifically,

$$\sigma_i(A) = \begin{cases} \sigma_i(X) - 1, & \text{if } i \in a, \\ 0 & \text{if } i \in b. \end{cases} \tag{30}$$

Clearly $A$ has $r$ numbers of nonzero distinct singular values. Denote them by $\nu_1(A) > \ldots > \nu_r(A)$. The index sets $a_1, \ldots, a_r$ that depending on $X$ in (22) also provides a partition of $(\sigma_i(A))_{i\in a}$, i.e.,

$$\sigma_i(A) = \nu_l(A), \quad \forall\, i \in a_l,\ \forall\, l = 1, \ldots, r.$$

For $l = 1, \ldots, r$, denote $\Omega_{al} : \mathbb{R}^{m\times n} \times \mathbb{R}^{m\times n} \to \mathbb{R}^{a_l \times a_l}$ as

$$\begin{aligned} \Omega_{a_l}(A, H) :=\ & (S(\widetilde{H}_1))_{a_l}^T(\mathrm{Diag}(\sigma(A)) - \nu_l(A)I_m)^\dagger(S(\widetilde{H}_1))_{a_l} + (2\nu_l(A))^{-1}\widetilde{H}_{a_lc}\widetilde{H}_{a_lc}^T \\ & + (T(\widetilde{H}_1))_{a_l}^T(-\mathrm{Diag}(\sigma(A)) - \nu_l(A)I_m)^\dagger(T(\widetilde{H}_1))_{a_l}, \quad (A, H) \in \mathbb{R}^{m\times n} \times \mathbb{R}^{m\times n}, \end{aligned}$$

where $I_m$ is the $m \times m$ identity matrix and $Z^\dagger$ denotes the Moore-Penrose pseudoinverse of a given matrix $Z$. Then for any $H \in \mathbb{R}^{m\times n}$, the conjugate of $\theta''(A; H, \cdot)$ is

$$\psi^*_{(A,H)}(B) := (\theta''(A; H, \cdot))^*(B) = 2\sum_{l=1}^{r} \mathrm{tr}\big(\Omega_{a_l}(A, H)\big) + 2\langle \mathrm{Diag}(\sigma_b(B)), U_b^T H A^\dagger H V_b\rangle, \tag{31}$$

where $\sigma_b(B) = (\sigma_i(B))_{i\in b}$ and $\widetilde{H} = [\widetilde{H}_1\ \widetilde{H}_2] = [U^T H V_1\ U^T H V_2]$.

In the following, we present several properties regarding the critical cone of $\theta$ and the directional derivative of $\mathrm{Prox}_\theta(\cdot)$.

**Proposition 4.1.** *Suppose that $X \in \mathbb{R}^{m\times n}$ has the singular value decomposition (20). Let the index sets $a, b, a_1, \ldots, a_l, b_1, b_2, b_3$ be given by (21) and (22). Given any $H \in \mathbb{R}^{m\times n}$, denote $\widetilde{H} = U^T H V$ for $(U, V) \in \mathcal{O}^{m\times n}(X)$. Denote $A = \mathrm{Prox}_\theta(X)$ and $B = \mathrm{Prox}_{\theta^*}(X)$. Then the following conclusions hold:*

*(i) $H \in \mathcal{C}_\theta(A, B)$ if and only if $H$ satisfies*

$$\widetilde{H} = \left(\begin{array}{c|ccc|c} \widetilde{H}_{aa} & & \widetilde{H}_{ab} & & \widetilde{H}_{ac} \\ \hline & \Pi_{\mathcal{S}_+^{|b_1|}}(\widetilde{H}_{b_1b_1}) & 0 & 0 & \\ \widetilde{H}_{ba} & 0 & 0_{b_2\times b_2} & 0 & 0 \\ & 0 & 0 & 0_{b_3\times b_3} & \end{array}\right). \tag{32}$$

*(ii) For any $D \in \mathbb{R}^{m\times n}$, $H = \mathrm{Prox}'_\theta(X; H+D)$ if and only if $H \in \mathcal{C}_\theta(A, B)$ and $\langle H, D\rangle = -\psi^*_{(A,H)}(B)$, where the function $\psi^*_{(A,H)}(\cdot)$ is given in (31).*

*Proof.* The result of part (i) can be obtained from [6, proposition 10]. Now we derive (iii). Suppose $H = \mathrm{Prox}'_\theta(C; H + D)$. Denote $\widetilde{H} = [\widetilde{H}_1, \widetilde{H}_2]$ with $\widetilde{H}_1 \in \mathbb{R}^{m\times m}$ and $\widetilde{H}_2 \in \mathbb{R}^{m\times(n-m)}$. Direct computations of $\psi^*_{(A,H)}(B)$ given in (31) show that

$$\begin{aligned}\psi^*_{(A,H)}(B) = & \sum_{1\leqslant l,t\leqslant r}\frac{2}{-\nu_t(A)-\nu_l(A)}\|(T(\widetilde{H}_1))_{a_la_t}\|^2+\sum_{1\leqslant l\leqslant r}\frac{4}{-\nu_l(A)}\|(T(\widetilde{H}_1))_{a_lb_1}\|^2\\ &+\sum_{\substack{1\leqslant l\leqslant r\\ 1\leqslant i-|a|-|b_1|\leqslant|b_2|}}\left(\frac{2(1-\sigma_i(B))}{-\nu_l(A)}\|(S(\widetilde{H}_1))_{a_li}\|^2+\frac{2(\sigma_i(B)+1)}{-\nu_l(A)}\|(T(\widetilde{H}_1))_{a_li}\|^2\right)\\ &+\sum_{\substack{1\leqslant l\leqslant r\\ 1\leqslant i-|a|-|b_1|-|b_2|\leqslant|b_3|}}\left(\frac{2}{-\nu_l(A)}\|(S(\widetilde{H}_1))_{a_li}\|^2+\frac{2}{-\nu_l(A)}\|(T(\widetilde{H}_1))_{a_li}\|^2\right)\\ &+\sum_{1\leqslant l\leqslant r}\frac{1}{-\nu_l(A)}\|(\widetilde{H}_c)_{a_lc}\|^2.\end{aligned}\tag{33}$$

Recall the formula of $\mathrm{Prox}'_\theta(X;\cdot)$ given in (25). We deduce that

$$\widetilde{H}=\left(\begin{array}{c|ccc|c}\widetilde{H}_{aa} & & \widetilde{H}_{ab} & & \widetilde{H}_{ac}\\ \hline & \widetilde{H}_{b_1b_1} & 0 & 0 & \\ \widetilde{H}_{ba} & 0 & 0_{b_2\times b_2} & 0 & 0\\ & 0 & 0 & 0_{b_3\times b_3} & \end{array}\right)$$

and

$$\begin{cases}\widetilde{D}_{a_la_t}=\dfrac{2}{\nu_l(X)+\nu_t(X)-2}(T(\widetilde{H}_1))_{a_la_t},\ 1\leqslant l,t\leqslant r,\\ (\widetilde{D}_{ab})_{ij}=\dfrac{1}{\sigma_i(X)-1}(\widetilde{H}_{ab})_{ij}-\dfrac{\sigma_{j+|a|}(X)}{\sigma_i(X)-1}(\widetilde{H}_{ab})_{ji},\ i=1,2,\dots|a|,\ j=1,2,\dots,|b|,\\ (\widetilde{D}_{ba})_{ji}=\dfrac{1}{\sigma_i(X)-1}(\widetilde{H}_{ab})_{ji}-\dfrac{\sigma_{j+|a|}(X)}{\sigma_i(X)-1}(\widetilde{H}_{ab})_{ij},\ i=1,2,\dots|a|,\ j=1,2,\dots,|b|,\\ (\widetilde{D}_{ac})_{ij}=\dfrac{1}{\sigma_i(X)-1}(\widetilde{H}_{ac})_{ij},\ i=1,2,\dots,|a|,\ j=1,2,\dots,n-m,\\ \mathbb{S}^{|b_1|}_+\ni\widetilde{H}_{b_1b_1}=S(\widetilde{H}_{b_1b_1})\perp S(\widetilde{D}_{b_1b_1})\in\mathbb{S}^{|b_1|}_-,\end{cases}$$

where $\widetilde{D}=U^TDV$. Consequently, it follows from part (i) of this proposition that $H\in\mathcal{C}_\theta(A,B)$. Moreover, we have

$$\begin{aligned}\langle D,H\rangle &=\langle\tilde{d}_{aa},\widetilde{H}_{aa}\rangle+\langle\tilde{d}_{ab},\widetilde{H}_{ab}\rangle+\langle\tilde{d}_{ba},\widetilde{H}_{ba}\rangle+\langle\tilde{d}_{a2},\widetilde{H}_{a2}\rangle\\ &=\sum_{1\leqslant t,l\leqslant r}\frac{2}{\nu_l(X)+\nu_t(X)-2}\|(T(\widetilde{H}_1))_{a_la_t}\|^2+\sum_{1\leqslant l\leqslant r}\frac{4}{\nu_l(X)-1}\|(T(\widetilde{H}_1))_{a_lb_1}\|^2\\ &+\sum_{\substack{1\leqslant l\leqslant r\\ 1\leqslant i-|a|-|b_1|\leqslant|b_2|}}\left(\frac{2(1-\sigma_i(X))}{\nu_l(X)-1}\|(S(\widetilde{H}_1))_{a_li}\|^2+\frac{2(\sigma_i(X)+1)}{\nu_l(X)-1}\|(T(\widetilde{H}_1))_{a_li}\|^2\right)\\ &+\sum_{\substack{1\leqslant l\leqslant r\\ 1\leqslant i-|a|-|b_1|-|b_2|\leqslant|b_3|}}\left(\frac{2}{\nu_l(X)-1}\|(S(\widetilde{H}_1))_{a_li}\|^2+\frac{2}{\nu_l(X)-1}\|(T(\widetilde{H}_1))_{a_li}\|^2\right)\\ &\quad+\sum_{1\leqslant l\leqslant r}\frac{1}{\nu_l(X)-1}\|(\widetilde{H}_2)_{a_lc}\|^2.\end{aligned}$$

Note from (30) that $\nu_l(A) = \nu_l(X) - 1$ for any $l = 1, \ldots, r$. Hence, $\langle D, H \rangle = -\psi^*_{(A,H)}(B)$ by (33) and the above equation. The converse of this statement can be established by reversing the above arguments. $\square$

**Proposition 4.2.** *Suppose that $X \in \mathbb{R}^{m\times n}$ has the singular value decomposition (20). Let the index sets $a, b, a_1, \ldots, a_l, b_1, b_2, b_3$ be given by (21) and (22). Given any $H \in \mathbb{R}^{m\times n}$, denote $\widetilde{H} = U^T HV$ for $(U, V) \in \mathcal{O}^{m\times n}(X)$. Denote $A = \mathrm{Prox}_\theta(X)$ and $B = \mathrm{Prox}_{\theta^*}(X)$. Then the following conclusions hold:*

*(i) $H \in \mathcal{C}_{\theta^*}(B, A)$ if and only if $H$ satisfies $S(\widetilde{H}_{b_1b_1}) \in \mathbb{S}^{|b_1|}_-$ and*

$$\widetilde{H} = \left(\begin{array}{c|c|c|c|c} T(\widetilde{H}_{aa}) & \frac{1}{2}(\widetilde{H}_{ab_1} - \widetilde{H}^T_{b_1a}) & \widetilde{H}_{ab_2} & \widetilde{H}_{ab_3} & \\ \cline{1-4} \frac{1}{2}(\widetilde{H}_{b_1a} - \widetilde{H}^T_{ab_1}) & \widetilde{H}_{b_1b_1} & \widetilde{H}_{b_1b_2} & \widetilde{H}_{b_1b_3} & \\ \cline{1-4} \widetilde{H}_{b_2a} & \widetilde{H}_{b_2b_1} & \widetilde{H}_{b_2b_2} & \widetilde{H}_{b_2b_3} & \widetilde{H}_c \\ \cline{1-4} \widetilde{H}_{b_3a} & \widetilde{H}_{b_3b_1} & \widetilde{H}_{b_3b_2} & \widetilde{H}_{b_3b_3} & \end{array}\right), \tag{34}$$

*where the linear operators $S(\cdot)$ and $T(\cdot)$ are defined in (24).*

*(ii) $H \in (\mathcal{C}_\theta(A, B))^\circ$ if and only if $\phi^*_{(B,H)}(A) = 0$ and $H \in \mathcal{C}_{\theta^*}(B, A)$, where $\phi^*_{(B,H)}(\cdot)$ is the conjugate function of $\phi_{(B,H)}(\cdot) = (\theta^*)''(B; H, \cdot)$.*

*(iii) $H \in (\mathcal{C}_{\theta^*}(B, A))^\circ$ if and only if $\psi^*_{(A,H)}(B) = 0$ and $H \in \mathcal{C}_\theta(A, B)$.*

*Proof.* Part (i) follows from [6, Proposition 12]. To prove part (ii), we use a result from [6, Proposition 16] stating that $\phi^*_{(B,H)}(A) = 0$ if and only if $\psi^*_{(A,H)}(B) = 0$ for any $H \in \mathbb{R}^{m\times n}$, which is further equivalent to

$$\widetilde{H}_{aa} \in \mathbb{S}^{|a|}, \quad \widetilde{H}_{ab_1} = \widetilde{H}^T_{b_1a}, \quad \widetilde{H}_{ab_2} = \widetilde{H}^T_{b_2a} = 0, \quad \widetilde{H}_{ab_3} = \widetilde{H}^T_{b_3a} = 0, \quad \widetilde{H}_{ac} = 0.$$

Then by part (i) of this proposition, one can see that $\phi^*_{(B,H)}(A) = 0$ and $H \in \mathcal{C}_{\theta^*}(B, A)$ imply that

$$\widetilde{H}_{aa} = 0, \quad \widetilde{H}_{ab} = 0, \quad \widetilde{H}_{ac} = 0, \quad \widetilde{H}_{ba} = 0, \quad S(\widetilde{H}_{b_1b_1}) \in \mathbb{S}^{|b_1|}_-.$$

In view of Proposition 4.1, this is equivalent to have $H \in (\mathcal{C}_\theta(A, B))^\circ$. To prove part (iii), it suffices to note that either $H \in (\mathcal{C}_{\theta^*}(B, A))^\circ$ or $\psi^*_{(A,H)}(B) = 0$ with $H \in \mathcal{C}_\theta(A, B)$ is equivalent to

$$\widetilde{H} = \left(\begin{array}{c|c|c|c|c} S(\widetilde{H}_{aa}) & \frac{1}{2}(\widetilde{H}_{ab_1} + \widetilde{H}^T_{b_1a}) & 0 & 0 & \\ \cline{1-4} \frac{1}{2}(\widetilde{H}_{b_1a} + \widetilde{H}^T_{ab_1}) & \Pi_{\mathcal{S}^{|b_1|}_+}(\widetilde{H}_{b_1b_1}) & 0 & 0 & 0 \\ \cline{1-4} 0 & 0 & 0 & 0 & \\ \cline{1-4} 0 & 0 & 0 & 0 & \end{array}\right).$$

The proof of this proposition is completed. $\square$

## 5 The robust isolated calmness of the KKT solution mapping

The aim of this section is to show that the SOSC for the primal problem (2) (the dual problem (4)) is in fact equivalent to the SRCQ for the dual problem (4) (the primal problem (2)). Before proceeding, we mention that a variation of this result regarding the linear semidefinite programming has been studied in [20].

Following the notation in the previous section, we use $\theta$ to denote the nuclear norm function in $\mathbb{R}^{m\times n}$. Let $\Omega_P \subseteq \mathbb{X}$ and $\Omega_D \subseteq \mathbb{R}^e \times \mathbb{R}^d \times \mathbb{R}^{m\times n}$ be the optimal solution sets of the primal problem (2) and the dual problem (4), respectively, both being assumed to be non-empty. It follows from (14) that the KKT optimality condition of problem (2) is given by

$$\begin{cases} 0 \in \mathcal{F}^*\nabla h(\mathcal{F}X) + C + \partial\theta(X) + \mathcal{A}^*y, \\ y \in \mathcal{N}_{\mathcal{Q}}(\mathcal{A}X - b), \end{cases} \quad (X, y) \in \mathbb{R}^{m\times n} \times \mathbb{R}^e. \tag{35}$$

We write $\mathcal{M}_P(\overline{X}) \subseteq \mathbb{R}^e$ as the set of Lagrangian multipliers $\bar{y}$ associated with $\overline{X} \in \Omega_P$, i.e., $\bar{y} \in \mathcal{M}_P(\overline{X})$ if and only if $(\overline{X}, \bar{y})$ satisfies (35). Let $\mathcal{M}_D(\bar{y}, \bar{w}, \overline{S}) \subseteq \mathbb{R}^{m\times n}$ be the set of Lagrangian multipliers associated with $(\bar{y}, \bar{w}, \overline{S}) \in \Omega_D$ for problem (4), i.e., $\overline{X} \in \mathcal{M}_D(\bar{y}, \bar{w}, \overline{S})$ if and only if $(\overline{X}, \bar{y}, \bar{w}, \overline{S})$ solves the following KKT system:

$$\begin{cases} \mathcal{A}X - b \in \mathcal{N}_{\mathcal{Q}^\circ}(y), \ \ \mathcal{F}X \in \partial h^*(w), \\ X \in \partial\theta^*(S), \ 0 = \mathcal{A}^*y + \mathcal{F}^*w + S + C, \end{cases} \quad (y, w, S, X) \in \mathbb{R}^e \times \mathbb{R}^d \times \mathbb{R}^{m\times n} \times \mathbb{R}^{m\times n}. \tag{36}$$

Since $h$ is assumed to be essentially strictly convex, $h^*$ is essentially smooth [14, Theorem 26.3]. Thus, $\nabla h^*$ is locally Lipschitz continuous and directionally differentiable on $\mathrm{int}\,(\mathrm{dom}\, h^*)$. Moreover, $\partial h^*(w) = \varnothing$ whenever $w \notin \mathrm{int}\,(\mathrm{dom}\, h^*)$ [14, Theorem 26.1]. Therefore, if (36) admits at least one solution, this KKT optimality condition can be equivalently written as

$$\begin{cases} \mathcal{A}X - b \in \mathcal{N}_{\mathcal{Q}^\circ}(y), \ \ \mathcal{F}X \in \nabla h^*(w), \\ X \in \partial\theta^*(S), \ 0 = \mathcal{A}^*y + \mathcal{F}^*w + S + C, \end{cases} \quad (y, w, S, X) \in \mathbb{R}^e \times \mathbb{R}^d \times \mathbb{R}^{m\times n} \times \mathbb{R}^{m\times n}.$$

As in Section 3, we consider the canonical perturbation of problem (2) for the sake of subsequent sensitivity analysis:

$$\begin{array}{ll} \min\limits_{X} & h(\mathcal{F}X) + \langle C, X\rangle + \|X\|_* - \langle X, \delta_1\rangle \\ \text{s.t.} & \mathcal{A}X - b + \delta_2 \in \mathcal{Q}, \end{array}$$

where $\delta_1 \in \mathbb{R}^{m\times n}$ and $\delta_2 \in \mathbb{R}^e$ are perturbation parameters. For any given $(\delta_1, \delta_2) \in \mathbb{R}^{m\times n} \times \mathbb{R}^e$, the KKT optimality condition then takes the form of

$$\begin{cases} \delta_1 \in \mathcal{F}^*\nabla h(\mathcal{F}X) + C + \partial\theta(X) + \mathcal{A}^*y, \\ y \in \mathcal{N}_{\mathcal{Q}}(\mathcal{A}X - b + \delta_2), \end{cases} \quad (X, y) \in \mathbb{R}^{m\times n} \times \mathbb{R}^e. \tag{37}$$

Let $S_{\mathrm{KKT}} : \mathbb{R}^{m\times n} \times \mathbb{R}^e \to \mathbb{R}^{m\times n} \times \mathbb{R}^e$ be the following KKT solution mapping:

$$S_{\mathrm{KKT}}(\delta_1, \delta_2) := \{(x, y) \in \mathbb{R}^{m\times n} \times \mathbb{R}^e \mid (x, y) \text{ satisfies (37)}\}, \ (\delta_1, \delta_2) \in \mathbb{R}^{m\times n} \times \mathbb{R}^e. \tag{38}$$

The RCQ of problem (2) at a feasible solution $\overline{X} \in \mathbb{R}^{m\times n}$ is given by

$$\mathcal{A}\mathbb{R}^{m\times n} + \mathcal{T}_{\mathcal{Q}}(\mathcal{A}\overline{X} - b) = \mathbb{R}^e. \tag{39}$$

Let $(y, q) \in \mathbb{R}^e \times \mathbb{R}^e$ satisfy $y \in \mathcal{N}_{\mathcal{Q}}(q)$. We denote the critical cone of $\mathcal{Q}$ at $q$ for $y$ and the critical cone of $\mathcal{Q}^\circ$ at $y$ for $q$ as

$$\mathcal{C}_{\mathcal{Q}}(q, y) := \mathcal{T}_{\mathcal{Q}}(q) \cap y^\perp, \qquad \mathcal{C}_{\mathcal{Q}^\circ}(y, q) := \mathcal{T}_{\mathcal{Q}^\circ}(y) \cap q^\perp.$$

It is easy to verify that

$$(\mathcal{C}_{\mathcal{Q}}(q, y))^\circ = \mathcal{C}_{\mathcal{Q}^\circ}(y, q). \tag{40}$$

The following theorem, which is the main result of our paper, demonstrates the equivalence between the primal SOSC and the dual SRCQ.

**Theorem 5.1.** *Let $\overline{X} \in \mathbb{R}^{m\times n}$ be an optimal solution of problem (2) with $\mathcal{M}_P(\overline{X}) \neq \emptyset$. Let $\bar{y} \in \mathcal{M}_P(\overline{X})$. Denote $\overline{S} := -\mathcal{A}^*\bar{y} - \mathcal{F}^*\nabla h(\mathcal{F}\overline{X}) - C$. Then the following two statements are equivalent to each other:*

*(i) The SOSC holds at $\overline{X}$ for $\bar{y}$ with respect to the primal problem (2), i.e.,*

$$\langle \mathcal{F}H, \nabla^2 h(\mathcal{F}\overline{X})\mathcal{F}H\rangle - \psi^*_{(\overline{X},H)}(\overline{S}) > 0, \quad \forall\, H \in \mathcal{C}(\bar{x})\backslash\{0\}, \tag{41}$$

*where $\mathcal{C}(\overline{X}) := \mathcal{C}_{\mathcal{Q}}(\mathcal{A}\overline{X} - b, \bar{y}) \cap \mathcal{C}_\theta(\overline{X}, \overline{S})$.*

*(ii) The SRCQ holds at $\bar{y}$ for $\overline{X}$ with respect to the dual problem (4), i.e.,*

$$\mathcal{F}^*\mathbb{R}^d + \mathcal{A}^*\mathcal{C}_{\mathcal{Q}^\circ}(\bar{y}, \mathcal{A}\overline{X} - b) + \mathcal{C}_{\theta^*}(\overline{S}, \overline{X}) = \mathbb{R}^{m\times n} \tag{42}$$

*Proof.* Firstly, let us assume that the statement (i) holds. Denote

$$\mathbb{E} := \mathcal{F}^*\mathbb{R}^d + \mathcal{A}^*\mathcal{C}_{\mathcal{Q}^\circ}(\bar{y}, \mathcal{A}\overline{X} - b) + \mathcal{C}_{\theta^*}(\overline{S}, \overline{X}).$$

Suppose on the contrary that $\mathbb{E} \neq \mathbb{R}^{m\times n}$. Then $\mathrm{cl}(\mathbb{E}) \neq \mathbb{R}^{m\times n}$ [14, Theorem 6.3]. Hence, there exists $D \in \mathbb{R}^{m\times n}$ but $D \notin \mathrm{cl}(\mathbb{E})$. Note that $\mathrm{cl}(\mathbb{E})$ is a closed convex cone. Let

$$\overline{D} := D - \Pi_{\mathrm{cl}(\mathbb{E})}(D) = \Pi_{(\mathrm{cl}(\mathbb{E}))^\circ}(D) \neq 0.$$

Obviously, $\langle H, \overline{D}\rangle \leqslant 0$ for any $H \in \mathrm{cl}(\mathbb{E})$. This implies that

$$\mathcal{F}\overline{D} = 0, \quad \mathcal{A}\overline{D} \in (\mathcal{C}_{\mathcal{Q}^\circ}(\bar{y}, \mathcal{A}\overline{X} - b))^\circ, \quad \overline{D} \in (\mathcal{C}_{\theta^*}(\overline{S}, \overline{X}))^\circ.$$

Thus, it follows from (40) that $\mathcal{A}\overline{D} \in \mathcal{C}_{\mathcal{Q}}(\mathcal{A}\overline{X} - b, \bar{y})$. From Proposition 4.2, we also have that

$$\psi^*_{(\overline{X},\overline{D})}(\overline{S}) = 0, \quad \overline{D} \in \mathcal{C}_\theta(\overline{X}, \overline{S}).$$

Therefore, $\overline{D} \in \mathcal{C}(\overline{X})\backslash\{0\}$ and $\langle \mathcal{F}\overline{D}, \nabla^2 h(\mathcal{F}\overline{X})\mathcal{F}\overline{D}\rangle - \psi^*_{(\overline{X},\overline{D})}(\overline{S}) = 0$, which contradicts the assumed SOSC (41) at $\overline{X}$.

The reverse implication can be proved similarly. Suppose that the SOSC (41) fails to hold at $\overline{X}$ for $\bar{y}$. Then there exists $\overline{H} \in \mathcal{C}(\overline{X})\backslash\{0\}$ such that

$$\langle \mathcal{F}\overline{H}, \nabla^2 h(\mathcal{F}\overline{X})\mathcal{F}\overline{H}\rangle - \psi^*_{(\overline{X},\overline{H})}(\overline{S}) = 0.$$

Since $h$ is assumed to be essentially strictly convex, $\langle \mathcal{F}H, \nabla^2 h(\mathcal{F}\overline{X})\mathcal{F}H\rangle > 0$ for any $H \in \mathbb{R}^{m\times n}$ such that $\mathcal{F}H \neq 0$. It also follows from [6, Proposition 16] that $\psi^*_{(\overline{X},H)}(\overline{S}) \leqslant 0$ for any $H \in \mathbb{R}^{m\times n}$. Consequently,

$$\mathcal{F}\overline{H} = 0, \quad \psi^*_{(\overline{X},\overline{H})}(\overline{S}) = 0.$$

We have from $\overline{H} \in \mathcal{C}(\overline{X})\backslash\{0\}$ that

$$\mathcal{A}\overline{H} \in \mathcal{C}_{\mathcal{Q}}(\mathcal{A}\overline{X} - b, \bar{y}), \quad \overline{H} \in \mathcal{C}_{\theta}(\overline{X}, \overline{S}).$$

Hence, we deduce from (40) and Proposition 4.2 that

$$\overline{H} \in (\mathcal{A}^*\mathcal{C}_{\mathcal{Q}^\circ}(\bar{y}, \mathcal{A}\overline{X} - b))^\circ \cap (\mathcal{C}_{\theta^*}(\overline{S}, \overline{X}))^\circ = (\mathcal{A}^*\mathcal{C}_{\mathcal{Q}^\circ}(\bar{y}, \mathcal{A}\overline{X} - b) + \mathcal{C}_{\theta^*}(\overline{S}, \overline{X}))^\circ.$$

By the assumed SRCQ (42) at $\bar{y}$ for $\overline{X}$, there exist $\tilde{w} \in \mathbb{R}^d$ and $\widetilde{H} \in \mathcal{A}^*\mathcal{C}_{\mathcal{Q}^\circ}(\bar{y}, \mathcal{A}\overline{X} - b) + \mathcal{C}_{\theta^*}(\overline{S}, \overline{X})$ such that $\overline{H} = \mathcal{F}^*\tilde{w} + \widetilde{H}$. Then

$$\langle \overline{H}, \overline{H} \rangle = \langle \overline{H}, \mathcal{F}^*\tilde{w} + \widetilde{H} \rangle = \langle \overline{H}, \widetilde{H} \rangle \leqslant 0,$$

which implies $\overline{H} = 0$. This contradicts the previous assumption that $\overline{H} \neq 0$. The proof is thus completed. □

One can also establish an analogous result by swapping the roles of the primal and dual problems in Theorem 5.1.

**Theorem 5.2.** *Let $(\bar{y}, \bar{w}, \overline{S}) \in \mathbb{R}^e \times \mathbb{R}^d \times \mathbb{R}^{m \times n}$ be an optimal solution of problem (4) with $\mathcal{M}_D(\bar{y}, \bar{w}, \overline{S}) \neq \varnothing$. Let $\overline{X} \in \mathcal{M}_D(\bar{y}, \bar{w}, \overline{S})$. Then the following two statements are equivalent to each other:*

*(i) The SOSC holds at $(\bar{y}, \bar{w}, \overline{S})$ for $\overline{X}$ with respect to the dual problem (4), i.e.,*

$$\langle H_w, (\nabla h^*)'(\bar{w}; H_w) \rangle - \phi^*_{(\overline{S}, H_S)}(\overline{X}) > 0, \quad \forall\ (H_y, H_w, H_s) \in \mathcal{C}(\bar{y}, \bar{w}, \overline{S})\backslash\{0\}, \tag{43}$$

*where the critical cone $\mathcal{C}(\bar{y}, \bar{w}, \overline{S})$ is defined as*

$$\mathcal{C}(\bar{y}, \bar{w}, \overline{S}) := \left\{ (H_y, H_w, H_S) \in \mathbb{R}^e \times \mathbb{R}^d \times \mathbb{R}^{m \times n} \,\middle|\, \begin{array}{c} \mathcal{A}^*H_y + \mathcal{F}^*H_w + H_S = 0, \\ H_y \in \mathcal{C}_{\mathcal{Q}^\circ}(\bar{y}, \mathcal{A}\overline{X} - b),\ H_S \in \mathcal{C}_{\theta^*}(\overline{S}, \overline{X}) \end{array} \right\}.$$

*(ii) The SRCQ holds at $\overline{X}$ for $(\bar{y}, \bar{w}, \overline{S})$ with respect to the primal problem (2), i.e.,*

$$\begin{pmatrix} \mathcal{A} \\ \mathcal{I}_{\mathbb{R}^{m \times n}} \end{pmatrix} \mathbb{R}^{m \times n} + \begin{pmatrix} \mathcal{C}_{\mathcal{Q}}(\mathcal{A}\overline{X} - b, \bar{y}) \\ \mathcal{C}_{\theta}(\overline{X}, \overline{S}) \end{pmatrix} = \begin{pmatrix} \mathbb{R}^e \\ \mathbb{R}^{m \times n} \end{pmatrix}. \tag{44}$$

*Proof.* With the help of Proposition 4.2, one can establish the assertion in the same fashion as in Theorem 5.1. We omit the details here. □

Finally, by combining Proposition 3.3, Theorem 5.1 and Theorem 5.2, we are ready to provide several equivalent characterizations of the robust isolated calmness of the KKT solution mapping at the origin for the unique KKT point of problem (2).

**Theorem 5.3.** *Let $\overline{X} \in \mathbb{R}^{m \times n}$ be an optimal solution of problem (2) and $(\bar{y}, \bar{w}, \overline{S}) \in \mathbb{R}^e \times \mathbb{R}^d \times \mathbb{R}^{m \times n}$ be an optimal solution of problem (4). Assume that the RCQ (39) holds at $\overline{X}$. Then the following statements are equivalent to each other:*

*(i) The KKT solution mapping $S_{KKT}$ in (38) is robustly isolated calm at the origin for $(\overline{X}, \bar{y})$.*

*(ii) The SOSC (41) holds at $\overline{X}$ for $\bar{y}$ with respect to the primal problem (2) and the SRCQ (44) holds at $\overline{X}$ for $(\bar{y}, \bar{w}, \overline{S})$ with respect to the primal problem (2).*

*(iii) The SOSC (41) holds at* $\overline{X}$ *for* $\bar{y}$ *with respect to the primal problem (2) and the SOSC (43) holds at* $(\bar{y}, \bar{w}, \overline{S})$ *for* $\overline{X}$ *with respect to the dual problem (4).*

*(iv) The SRCQ (42) holds at* $\bar{y}$ *for* $\overline{X}$ *with respect to the dual problem (4) and the SRCQ (44) holds at* $\overline{X}$ *for* $(\bar{y}, \bar{w}, \overline{S})$ *with respect to the primal problem (2).*

*(v) The SOSC (43) holds at* $(\bar{y}, \bar{w}, \overline{S})$ *for* $\overline{X}$ *with respect to the dual problem (4) and the SRCQ (42) holds at* $\bar{y}$ *for* $\overline{X}$ *with respect to the dual problem (4).*

It is worth mentioning that in this paper, we focus on the characterizations of the robust isolated calmness of the KKT solution mapping for problem (2) when it admits a unique KKT point. It would be certainly interesting to know to what extent our results can be extended to the case when problem (2) admits non-unique solutions. We shall leave this as a future research topic.

# References

[1] J. F. Bonnans and A. Shapiro, *Perturbation Analysis of Optimization Problems*, Springer, New York, 2000.

[2] E. J. Candès and B. Recht, *Exact matrix completion via convex optimization*, Foundations of Computational Mathematics, 9 (2008), pp. 717–772.

[3] E. J. Candès and T. Tao, *The power of convex relaxation: near-optimal matrix completion*, IEEE Transactions on Information Theory, 56 (2009), pp. 2053–2080.

[4] F. H. Clarke, *Optimization and Nonsmooth Analysis*, Second edition, Classics in Applied Mathematics, 5, Society for Industrial and Applied Mathematics, Philadelphia, PA, 1990.

[5] C. Ding, *An Introduction to a Class of Matrix Optimization Problems*, Ph.D Thesis, Department of Mathematics, National University of Singapore, 2012.

[6] C. Ding, *Variational analysis of the Ky Fan k-norm*, Set-Valued and Variational Analysis, (2017), pp. 1–34, DOI:10.1007/s11228-016-0378-3.

[7] C. Ding, D. F. Sun, J. Sun, and K.-C. Toh, *Spectral operators of matrices*, arXiv:1401.2269, 2014.

[8] C. Ding, D. F. Sun, and L. W. Zhang, *Characterization of the robust isolated calmness for a class of conic programming problems*, SIAM Journal on Optimization, 27 (2017), pp. 67–90.

[9] A. L. Dontchev and R. T. Rockafellar, *Implicit Functions and Solution Mappings*, Springer, New York, 2009.

[10] D. R. Han, D. F. Sun, and L. W. Zhang, *Linear rate convergence of the alternating direction method of multipliers for convex composite programming*, arXiv:1508.02134, 2015; revised in 2016.

[11] Y. L. Liu and S. H. Pan, *Locally upper Lipschitz of the perturbed KKT system of Ky Fan k-norm matrix conic optimization problems*, arXiv:1509.00681, 2015.

[12] K. Löwner, *Über monotone matrixfunktionen*, Mathematische Zeitschrift, 38 (1934), pp. 177–216.

[13] B. Recht, M. Fazel, and P. A. Parrilo, *Guaranteed minimum rank solutions to linear matrix equations via nuclear norm minimization*, SIAM Review, 52 (2010), pp. 471–501.

[14] R. T. Rockafellar, *Convex Analysis*, Princeton University Press, 1970.

[15] R. T. Rockafellar, *Augmented Lagrangians and applications of the proximal point algorithm in convex programming*, Mathematics of Operations Research, 1 (1976), pp. 97–116.

[16] T. K. Pong, P. Tseng, S. W. Ji, and J. P. Ye, *Trace norm regularization: reformulations, algorithms, and multi-task learning*, SIAM Journal on Optimization, 20 (2010), pp. 3465–3489.

[17] G. A. Watson, *Characterization of the subdifferential of some matrix norms*, Linear Algebra and its Applications, 170 (1992), pp. 33–45.

[18] G. A. Watson, *On matrix approximation problems with Ky Fan k norms*, Numerical Algorithms, 5 (1993), pp. 263–272.

[19] M. Yuan, A. Ekici, Z. S. Lu, and R. D. C. Monteiro, *Dimension reduction and coefficient estimation in multivariate linear regression*, Journal of the Royal Statistical Society: Series B 69 (2007), pp. 329–346.

[20] X. Y. Zhao, D. F. Sun, and K.-C. Toh, *A Newton-CG augmented Lagrangian method for semidefinite programming*, SIAM Journal on Optimization, 20 (2010), pp. 1737–1765.

[21] Y. L. Zhang and L. W. Zhang, *On the upper Lipschitz property of the KKT mapping for nonlinear semidefinite optimization*, Operations Research Letters, 44 (2016), pp. 474–478.